\documentclass{amsart}
\usepackage[utf8]{inputenc}

\usepackage{amsmath,amsfonts, amssymb}
\usepackage{graphicx}
\usepackage{xcolor}
\usepackage{enumerate}
\usepackage{enumitem}
\usepackage{cite}
\usepackage{tabu}

\usepackage{stmaryrd} 

\usepackage{tikz}

\usepackage{fullpage}
\usepackage[active]{srcltx}

\usepackage{changepage}

\usepackage{chngcntr}

\usepackage{multibib}

\usepackage{url}
\usepackage{hyperref}
\hypersetup{
	colorlinks,
	linkcolor={red!50!black},
	citecolor={blue!50!black},
	urlcolor={blue!80!black}
}  

\theoremstyle{plain}
\newtheorem{theorem}{Theorem}[section]
\newtheorem*{theo:sutmanwithnonisolatinggutsishoriprime}{Theorem~\ref{theo:sutmanwithnonisolatinggutsishoriprime}}
\newtheorem*{thm:handlenumberofgutsishandlenumberofmfld}{Theorem~\ref{thm:handlenumberofgutsishandlenumberofmfld}}
\newtheorem*{thm:bothbounds}{Theorem~\ref{thm:bothbounds}}

\newtheorem{cor}[theorem]{Corollary}

\newtheorem{lemma}[theorem]{Lemma}
\newtheorem{claim}[theorem]{Claim}

\newtheorem{prob}[theorem]{Problem}

\theoremstyle{definition}

\newcommand{\comment}[1]{}

\newcommand{\bdry}{\ensuremath{\partial}}

\newcommand{\nbhd}{\ensuremath{N}}

\newcommand{\MN}{{\ensuremath{\rm{MN}}}}

\newcommand{\cut}{\ensuremath{\backslash}}

\definecolor{amaranth}{rgb}{0.9, 0.17, 0.31} 
\definecolor{carrotorange}{rgb}{0.93, 0.57, 0.13} 
\definecolor{citrine}{rgb}{0.89, 0.82, 0.04} 
\definecolor{dartmouthgreen}{rgb}{0.05, 0.5, 0.06} 
\definecolor{ballblue}{rgb}{0.13, 0.67, 0.8} 
\definecolor{ceruleanblue}{rgb}{0.16, 0.32, 0.75} 
\definecolor{amethyst}{rgb}{0.6, 0.4, 0.8} 
\definecolor{amber}{rgb}{1.0, 0.75, 0.0} 
\definecolor{burlywood}{rgb}{0.87, 0.72, 0.53} 
\newcommand{\ken}[1]{{\color{amaranth} #1}}
\newcommand{\fab}[1]{{\color{ceruleanblue} #1}}

\title{Handle number is not always realized by a minimal genus Seifert surface.}

\author{Kenneth L. Baker}

\address{Department of Mathematics, University of Miami, Coral Gables, FL 33146, USA}
\email{k.baker@math.miami.edu}

\author{Fabiola Manjarrez-Guti\'errez}
\address{Instituto de Matem\'aticas, Universidad Nacional Aut\'onoma de Mexico, Cuernavaca, Mor., MEXICO}
\email{fabiola.manjarrez@im.unam.mx}


\begin{document}

\begin{abstract}
We construct genus one knots whose handle number is only realized by Seifert surfaces of non-minimal genus.  These are counterexamples to the conjecture that the Seifert genus of a knot is its Morse-Novikov genus. As the Morse-Novikov genus may be greater than the Seifert genus, we define the {\em genus $g$ Morse-Novikov number} $\MN_g(L)$ as the minimum handle number among Seifert surfaces for $L$ of genus $g$. Since, as we further show, the Morse-Novikov genus and the minimal genus Morse-Novikov number are additive under connected sum of knots, it then follows that there exists examples for which the discrepancies between Seifert genus and Morse-Novikov genus and between the Morse-Novikov number and the minimal genus Morse-Novikov number can be made arbitrarily large.
\end{abstract}

\maketitle



\section{Introduction}
For an oriented link $L$ in $S^3$, its {\em Morse-Novikov number} $\MN(L)$ is the handle number $h(M, \gamma_+, \xi)$ of the sutured manifold $(M, \gamma_+)$ where $M$ is the exterior of $L$, $R_+(\gamma_+)= \bdry M$ and $\xi \in H_2(M,\bdry M)$ is its Seifert class.
The handle number $h(F)$ of a surface $F$ representing $\xi$ is the handle number of $(M_F, \gamma_F)$, the sutured manifold resulting from decomposing $(M,\gamma_+)$ along $F$, and  $h(M,\gamma_+,\xi)$ is the minimum among these.
For details and background about handle numbers of links and sutured manifolds, see \cite{BMG-MorseNovikovHandleNumbers} and the earlier works of Goda \cite{goda-handlenumber, goda-murasugisum}\footnote{Goda's count of handle number for links and Seifert surfaces is half of our count.}.


The link's {\em Morse-Novikov genus} $g_{\MN}(L)$ is the minimum genus $g(F)$ among Seifert surfaces $F$ of $L$ realizing $h(M,\gamma_+, \xi)$.  The Morse-Novikov genus is always realized by an incompressible Seifert surface \cite{BakerAdditivityofMorseNovikov}.  Indeed, in the examples heretofore known, the Morse-Novikov genus is realized by some minimal genus Seifert surface (eg.\ \cite{goda-handlenumber,goda-murasugisum,freegenusone, BMG-nearlyfibered}). However, as Goda demonstrates,  the Morse-Novikov number is not necessarily realized by every minimal genus Seifert surface of a link \cite{goda-handlenumber}.
Regardless, it seemed not unreasonable to suspect that the Morse-Novikov genus is actually realized by a minimal genus Seifert surface.  Indeed in \cite{BMG-MorseNovikovHandleNumbers} we conjectured that $g_{\MN}(L)= g(L)$; see also \cite[Question 1.4]{BakerAdditivityofMorseNovikov}. However, in this note we construct counterexamples.  While we use Whitehead doubles to do a variation of trick of Lyon \cite{lyon} to facilitate our construction,  hyperbolic counterexamples should exist too, see Problem~\ref{prob:hypexamples}.

\begin{theorem}\label{thm:main}
    There exist genus one knots with Morse-Novikov genus at least two.  
\end{theorem}

More specifically, our examples have a unique genus one Seifert surface of handle number $6$ and an incompressible genus $2$ Seifert surface of handle number at most $4$.

Seeing that the Morse-Novikov genus may be greater than the Seifert genus prompts the following definition.   Given an oriented link $L$ and an integer $g \geq g(L)$, the {\em genus $g$ Morse-Novikov number} $\MN_g(L)$ is the minimum handle number among Seifert surfaces for $L$ of genus $g$.   In Lemma~\ref{lem:MNg-structure}, we show how $\MN_g(L)$ behaves as a function of $g$.

Thereafter we prove that Morse-Novikov genus and the genus $g$ Morse-Novikov numbers are additive under connected sum, Theorems~\ref{thm:MNgenusisadditive} and Theorem~\ref{thm:additivityofGenusgMorseNovikov}, in order to show there exists examples for which the discrepancies between Seifert genus and Morse-Novikov genus and between the Morse-Novikov number and the minimal genus Morse-Novikov number can be made arbitrarily large.

\begin{theorem}
\label{thm:arblargediscrepancy}
    For each positive integer $n$, there exists a knot $K_n$ with $g_{\MN}(K_n)-g(K_n)\geq n$ and $\MN_{g(K_n)}(K_n) - \MN(K_n) \geq 2n$.
\end{theorem}


    %

\subsection{Acknowledgements}

KLB was partially supported by the Simons Foundation (gifts \#523883 and \#962034  to Kenneth L.\ Baker). 
He also thanks the University of Pisa for their hospitality where a portion of this work was done.

FM-G was partially supported by Grant UNAM-PAPIIT IN113323.

\section{A basic construction}

\begin{figure}
    \centering
    \includegraphics[width=0.9\linewidth]{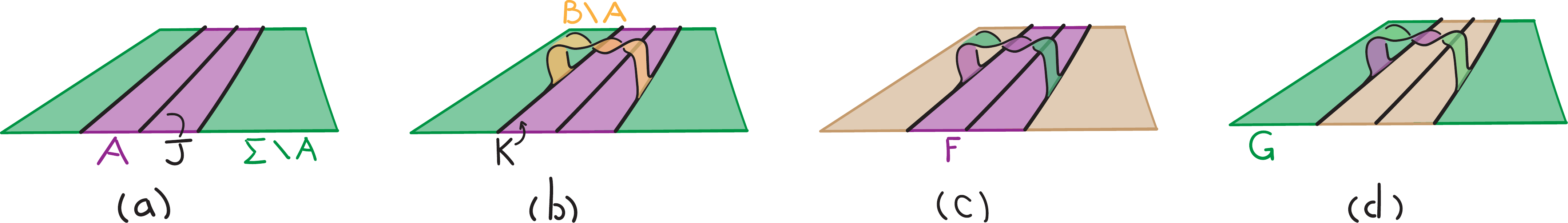}
    \caption{(a) A portion of $\Sigma$ near $J$ is shown with the annular neighborhood $A$ of $J$. (b) The portion $B\cut A$ of the Hopf band $B$ is highlighted and the knot $K$ is indicated. (c) \& (d) Attaching $B \cut A$ to $A$ or $\Sigma \cut A$ gives the Seifert surfaces $F$ and $G$.}
    \label{fig:surfacesABFGSimga}
\end{figure}

Let $J$ be a knot in a Heegaard surface $\Sigma$ and let $K$ be a $\Sigma$-framed Whitehead double of $J$. A genus one Seifert surface $F$ for $K$ is obtained from plumbing a Hopf band $B$ onto a closed annular neighborhood $A$ of $J$ in $\Sigma$.  Another Seifert surface for $K$ is $G = (\Sigma \cup B) \cut int(A)$ which has genus $g(G)=g(\Sigma)$.  See Figure~\ref{fig:surfacesABFGSimga}.
Observe that $F$ and $G$ may be isotoped rel-$K$ to have disjoint interiors so that the surface $F \cup G$ is a Heegaard stabilization of $\Sigma$.
(Using the orientations inherited from $\Sigma$, $F$ and $G$ induce opposite orientations on $K$.)
Do such an isotopy, and we may assume $(F \cup G) \cut  \nbhd(J) = \Sigma \cut \nbhd(J)$ so that $F \subset \nbhd(J)$ and $G \cap \nbhd(J)$ is a thrice punctured sphere.

\begin{theorem}
    Let $\sigma$ be the slope on $J$ induced by $\Sigma$.  
    Suppose that $J(\sigma)$, the result of $\sigma$-framed surgery on $J$, has Heegaard genus $HG(J(\sigma))$ equal to $g(\Sigma)+1$.  
    Further assume that the exterior of $J$ has no essential annuli with boundary slope $\sigma$.

    Then the Seifert surface $F$ is the unique genus $1$ Seifert surface for $K$, and $h(F)= 2g(\Sigma)+2$.
    However, the Seifert surface $G$ of genus $g(\Sigma)$ has $h(G)\leq 2g(\Sigma)$.  Consequently, $h(K) \leq h(G) < h(F)$.
\end{theorem}

\begin{proof}
    Suppose $K$ has a genus one Seifert surface that is not isotopic rel-$K$ to $F$.   Then there is one, say $F'$ that has interior disjoint from $F$ by Kakimizu \cite{KakimizuComplex}.  Since a component of the Whitehead link has a unique genus one Seifert surface disjoint from the other component, $F'$ must intersect the satellite torus $\bdry \nbhd(J)$.   In particular $F' \cut \nbhd(J)$ is an essential annulus in the exterior of $J$ with boundary slope $\sigma$, contrary to assumption.  Hence $F$ is the unique genus one Seifert surface for $K$.

    Since $F$ is obtained from $A$ by plumbing a Hopf band $B$, $h(F)=h(A)$ by \cite[Theorem 2]{goda-murasugisum} since the Hopf link fibers with fiber $B$. 
    The exterior of $A$ is a sutured manifold $(M_A,\gamma_A)$  which may be viewed as the exterior of $J$ with sutures that are a pair of annuli of slope $\sigma$. See Figure~\ref{fig:splittingforA}(a) which also shows the surface $\Sigma_A = \Sigma \cap M_A$. By definition, $h(A)=h(M_A,\gamma_A)$. 
    Then we have 
    \[h(F)=h(A)=h(M_A,\gamma_A) \geq 2 HG(J(\sigma)) = 2g(\Sigma)+2\]
    The inequality follows from \cite[Lemma 3.5]{BMG-nearlyfibered} since we chose $J$ so that $HG(J(\sigma))=g(\Sigma)+1$.

    To promote the inequality to an equality, one finds a Heegaard surface $\Sigma_A'$ that splits $(M_A,\gamma_A)$ into two compression bodies, each with handle number $g(\Sigma)+1$, as follows:  Take a meridional disk $D$ of $J$ so that $D \cut \nbhd(A)$ is an annulus that $\Sigma_A = \Sigma \cut \nbhd(A) \subset M_A$ intersects  in two spanning arcs.  Attach a handle to each side of $\Sigma_A$  along disjoint push-offs of the two arcs $\bdry D \cut \Sigma$ to form the surface $\Sigma_A'$ as in Figure~\ref{fig:splittingforA}(b).  
    Neglecting the sutured structure, since $\Sigma_A$ divides $M_A$ into two handlebodies of genus $g(\Sigma)$, we see that $\Sigma_A'$ divides $M_A$ into two handlebodies of genus $g(\Sigma)+2$.  (Each handlebody bounded by $\Sigma_A'$ is obtained as a handlebody bounded by $\Sigma$ with one handle attached outside and a boundary parallel arc drilled out from inside.)   However, with the sutured structures on these handlebodies arising by decomposing $(M_A,\gamma_A)$ along $\Sigma_A'$, each of the two disks $(D \cut \nbhd(A)) \cut \Sigma_A$ induce a product disk that runs once along the one of the attached handles (and crosses the corresponding annular component of $R(\gamma_A)$ once) as in Figure~\ref{fig:splittingforA}(c) \& (d). Thus each of these product decompositions yields a sutured manifold that is a handlebody of genus $g(\Sigma)+1$ for which the suture is a trivial curve bounding a disk, a ``spot''. See Figure~\ref{fig:splittingforA}(e). Hence these sutured manifolds are each compression bodies with handle number $g(\Sigma)+1$.  
    It then follows that $h(F)=2g(\Sigma)+2$.

    \begin{figure}
        \centering
        \includegraphics[width=0.9\linewidth]{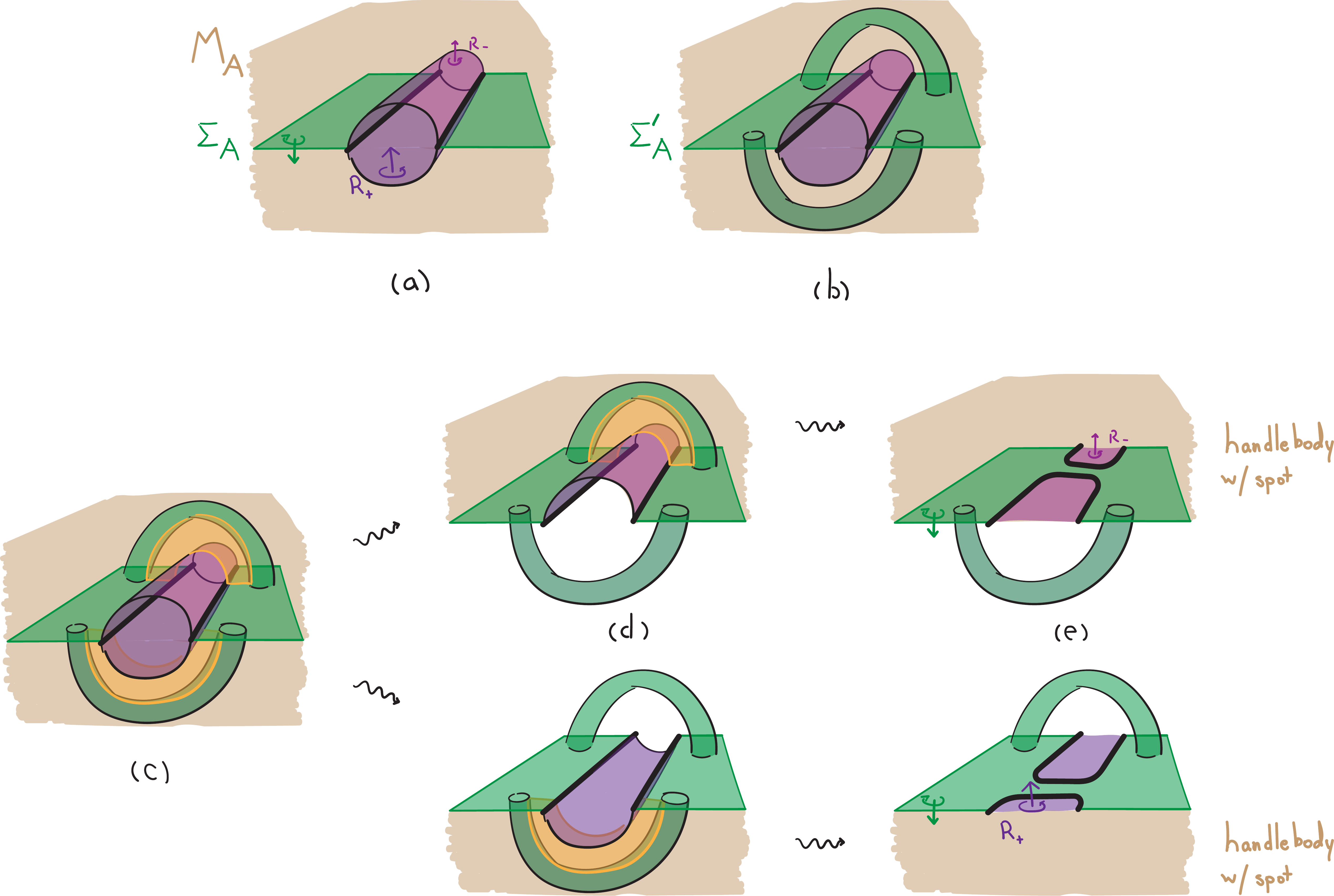}
        \caption{(a) The sutured manifold exterior $(M_A, \gamma_A)$ of the annulus $A$ is shown along with the surface $\Sigma_A$. (b) $\Sigma_A$ is tubed twice to form a Heegaard surface $\Sigma_A'$ for $(M_A, \gamma_A)$.  (c) \& (d) Two product disks for the decomposition of $M_A$ along $\Sigma_A'$. (e) The result of further decomposing along these product disks is a pair of simple compression bodies:  handlebodies with one ``spot''.
        }
        \label{fig:splittingforA}
    \end{figure}


    Now let $(M_G,\gamma_G)$ be the sutured manifold exterior of the Seifert surface $G = (\Sigma \cut A) \cup (B \cut A)$. Fixing $\Sigma \cut A$, isotop the band $B \cut A$ so that its core arc $b$ lies in $A \subset \Sigma$ as a spanning arc of $A$, as in the passage from Figure~\ref{fig:MGConstruction}(a) to (b).
    Choose a disk $D$ in $\Sigma$ that disjoint from $B$ and intersects $A$ as shown in upper portions of Figures~\ref{fig:MGConstruction}(c) \& (d).  Having $D$ overhang $A$ so that $D \cut A$ is a pair of bigons will be useful later in the construction.  Let us now regard $G$ as $\Sigma \cut (A \cup D)$ with the newly isotoped band $B \cut A$ attached.

    Consider a product neighborhood $\Sigma \times [-1,1]$ of $\Sigma = \Sigma \times \{0\}$.  
    Since $(\Sigma \cut A) \cup b \subset \Sigma$, the handlebody $\bar{\nbhd}(G)$ that is a closed regular neighborhood of $G$  lies in the interior of the product and may be viewed as the submanifold $(\Sigma\cut D) \times [-\epsilon, \epsilon]$ as in Figure~\ref{fig:MGConstruction}(e). 
    (Visualize this by fattening the band $B \cut A$ as done in the lower portions of Figures~\ref{fig:MGConstruction}(c) \& (d).)  
    Then the sutures $\gamma_G$ on $M_G = S^3 \cut \nbhd(G)$ may be regarded as the four arcs $((\bdry A) \cut D) \times \{\pm \epsilon\}$ joined by four spanning arcs of the annulus $\bdry D \times [-\epsilon, \epsilon]$ with a quarter twist.   Refer again to Figures~\ref{fig:MGConstruction}(c), (d), \& (e) to see how the arrangement of these sutures is obtained.
    
    \begin{figure}
        \centering
            \includegraphics[width=0.9\linewidth]{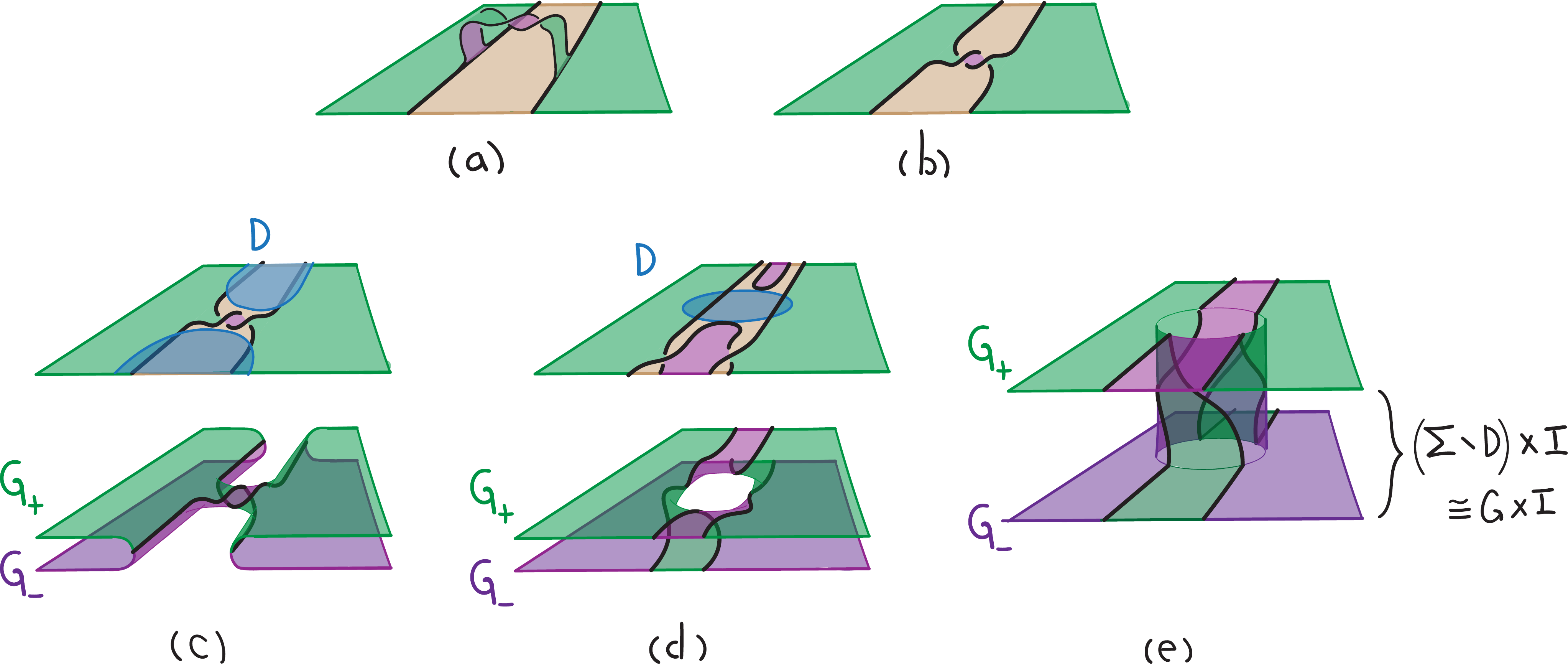}
        \caption{(a) The surface $G$ as a banding of $\Sigma \cut A$.  (b) The surface $G$ after isotoping the core arc of the band $B \cut A$ into $A$. (c)\& (d) A disk $D \subset \Sigma$ disjoint from $B\cut A$ (above) and the boundary of a corresponding handlebody thickening $\bar{\nbhd}(G)$ of $G$ (below). (e) A view of the boundary of the sutured manifold $(M_G, \gamma_G)$ complementary to $G$ in $(\Sigma \cut D) \times [-\epsilon, \epsilon]$.}
    \label{fig:MGConstruction}
    \end{figure}

    Form the Heegaard surface $S$ for $(M_G, \gamma_G)$ as follows, referring to Figure~\ref{fig:MGHeegaardSurface}.  Push copies $S_\pm$ of $(\Sigma \cut A) \times \{\pm \epsilon\}$ rel-$\bdry$ into $\Sigma \times [\epsilon,1]$ and $\Sigma \times [-1,-\epsilon]$ according to the sign $\pm$.  These surfaces will cross each of the disks $D \times \{\pm \epsilon\}$ in two arcs.  Those four arcs together with the four spanning arcs of $\gamma_G$ in the annulus  $\bdry D \times [-\epsilon, \epsilon]$  bound a disk $S_0$ in $D \times [-\epsilon, \epsilon]$.   This disk $S_0$ joins the two copies $S_+$ and $S_-$ of $\Sigma \cut A$ to form a surface $S$ whose boundary is the core curve of $\gamma_G$ by construction.  

    \begin{figure}
        \centering
            \includegraphics[width=0.9\linewidth]{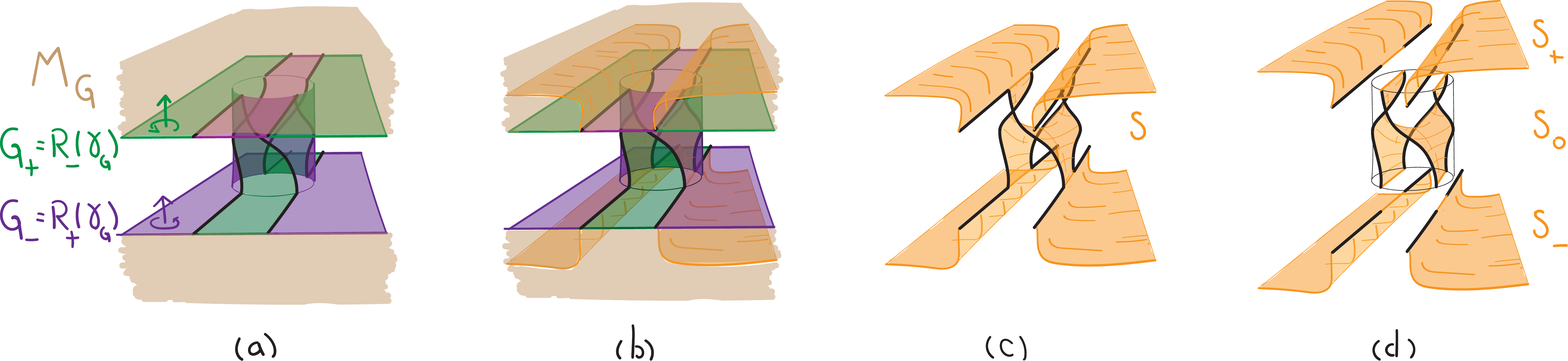}
        \caption{(a) The sutured manifold $(M_G, \gamma_G)$.  (b) The Heegaard surface $S$ as it sits inside $M_G$. (c) The Heegaard surface $S$ on its own. (d) An `exploded' view of $S$ for the pieces $S_+$, $S_0$, and $S_-$ of its construction.}
    \label{fig:MGHeegaardSurface}
    \end{figure}
    
    \begin{figure}
        \centering
        \includegraphics[width=0.9\linewidth]{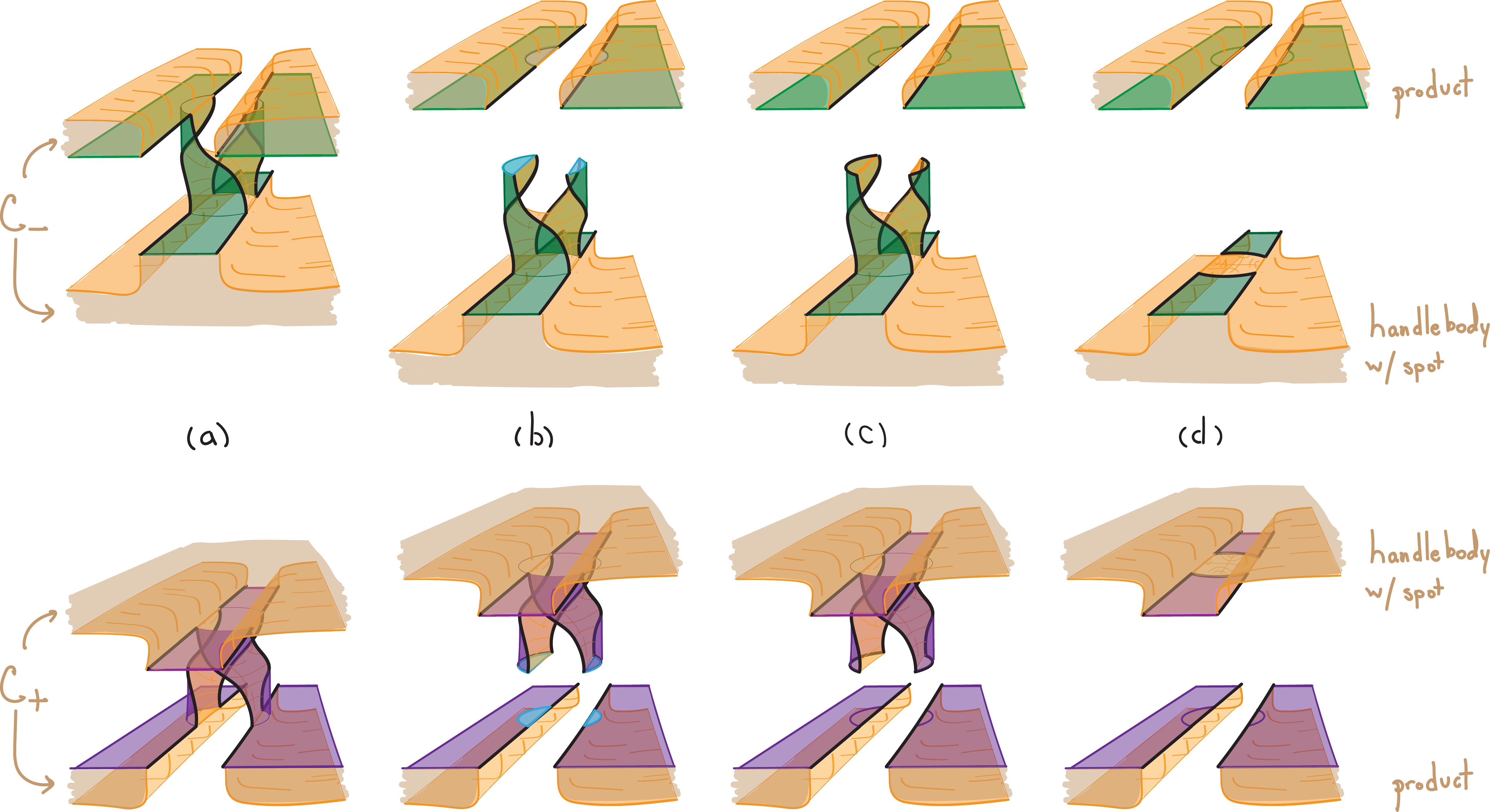}
        \caption{(a) The sutured manifolds $C_\pm$ resulting from decomposing $M_G$ along $S$. (b) A pair of product disks in each $C_\pm$. (c) The result of decomposing $C_\pm$ along these product disks. (d) An isotopy of the decompositions into more standard pieces. }
    \label{fig:MGDecomposition}
    \end{figure}

    The surface $S$ divides $(M_G, \gamma_G)$ into two sutured manifolds $(C_+, \gamma_{C_+})$ and $(C_-, \gamma_{C_-})$ shown separately in Figure~\ref{fig:MGDecomposition}(a).  The rest of Figure~\ref{fig:MGDecomposition} gives a decomposition of these sutured manifolds from which we may see that they are indeed compression bodies and determine their handle numbers.  The ensuing discussion of this below is phrased for $C_-$ but similarly applies to $C_+$.
    
    Here $C_-$ is the union of the handlebody below $S_- \cup A \times \{-\epsilon\}$ (which is the union of the handlebody below $\Sigma \times \{-1\}$ and the solid torus $A \times [-1, -\epsilon]$) joined to the product region between $(\Sigma \cut A)\times \{+\epsilon\}$ and $S_+$ by a component of $(D \times [-\epsilon, +\epsilon]) \cut S_0$. See Figure~\ref{fig:MGDecomposition}(b).  The suture $\gamma_{C_-}$ is determined by $\bdry S$ so that $S=R_+(\gamma_{C_-})$.

   Observe that $C_-$ intersects $D \times \{+\epsilon\}$ in two product disks, the disks shown in blue in  Figure~\ref{fig:MGDecomposition}(b).  Decomposing along these two product disks leaves a product sutured manifold between $(\Sigma \cut A) \times \{+\epsilon\}$ and $S_+$ and a compression body that is homeomorphic to the handlebody below $\Sigma \times \{-1\}$ in which the negative boundary is a disk that contains $(A\cut D)\times \{-\epsilon\}$,  ie.\ a handlebody with a `spot'. (The immediate result of the product disk decomposition is given in Figure~\ref{fig:MGDecomposition}(c) while Figure~\ref{fig:MGDecomposition}(d) shows the result after an isotopy.) The former has handle number $0$ while the latter has handle number equal to the genus of $\Sigma$.  Hence, by \cite[Lemma 2.4]{goda-murasugisum}, $C_-$ is a compression body and its handle number is the genus of $\Sigma$.
    
    

   The situation for the sutured manifold $(C_+,\gamma_{C_+})$ is analogous.  Hence $S$ is a Heegaard surface for $(M_G, \gamma_G)$ giving a splitting with handle number $2g(\Sigma)$.
\end{proof}

\section{Knots such as $J$ do exist. }

For integers $p,q,r,m,n$ such that  $p$ and $q$ are coprime, $m$ and $n$ are coprime, and a few other conditions are satisfied, Kang produces knots $K(p,q,r,m,n)$ in a genus $2$ Heegaard surface $\Sigma$ for $S^3$ such that $\Sigma$-framed surgery yields a Seifert fibered space over $S^2$ with four exceptional fibers of orders $p,q,m,n$, \cite[Theorem 3.6]{Kang}.  

\begin{lemma}
    Let $J$ be one of Kang's knots with parameters as above, and let $\sigma$ be its framing by $\Sigma$.   
    Then $HG(J(\sigma))=3$ and $\sigma$ is not an annular boundary slope for $J$.
\end{lemma}

\begin{proof}
The manifold $J(\sigma)$ is a Seifert fibered space over  $S^2$ with four exceptional fibers of orders $p,q,m,n$.  By the coprime conditions, at most two of these orders may be $2$.   Thus $J(\sigma)$ has Heegaard genus 3 by \cite{Morimoto}.  

If $\sigma$ were an annular boundary slope then $J$ would be a cabled knot with cabling slope $\sigma$.  Hence $J(\sigma)$ would be a reducible manifold with a lens space summand.  Yet that is contrary to $J(\sigma)$ being a Seifert fibered space over $S^2$ with four exceptional fibers.
\end{proof}

\begin{cor}\label{cor:Kangdoubles}
Let $K$ be a $\Sigma$-framed Whitehead double of one of Kang's knots as above.
Then $g(K)=1$, $\MN_1(K) = 6$, $g_\MN(K) \geq 2$, and $2 \leq \MN(K) \leq \MN_2(K) \leq 4$.  \qed



\end{cor}



\section{
Additivity of Morse-Novikov genus and minimal genus Morse-Novikov number}


\begin{lemma}[{\cite[Theorem 1 for $n=1$]{goda-murasugisum}}]\label{lem:handlenumberofboundarysumSeifertsurfaces}
    Let $S_a$ and $S_b$ be Seifert surfaces for knots $K_a$ and $K_b$.  Then $S_a \natural S_b$, the boundary sum of Seifert surfaces $S_i$ for $K_i$, is a Seifert surface for $K_a \# K_b$ and $h(S_a\natural S_b) = h(S_a) + h(S_b)$.  
    \qed
\end{lemma}


\begin{theorem}\label{thm:MNgenusisadditive}
    Morse-Novikov genus is additive under connected sum.

    Given knots $K_a$ and $K_b$,
    \[g_\MN(K_a \# K_b) = g_\MN(K_a) + g_\MN(K_b).\]

    More specifically,
    if a Seifert surface $S$ for $K=K_a \# K_b$ realizes 
     $g_{\MN}(K)$, then $S = S_a \natural S_b$ is the boundary sum of Seifert surfaces $S_i$ for $K_i$ that realize $g_{\MN}(K_i)$ for each $i=a,b$.
\end{theorem}

\begin{proof}
    Suppose a Seifert surface $S$ for $K=K_a \# K_b$ realizes $g_{\MN}(K)$.
    Then $S$ must be an incompressible Seifert surface \cite{BakerAdditivityofMorseNovikov}.
    Therefore the summing sphere realizing $K$ as $K_a \# K_b$ may be taken to intersect $S$ in just an arc realizing $S$ as $S_a \natural S_b$, the boundary sum of Seifert surfaces $S_i$ for $K_i$. Then $h(S)= h(S_a)+h(S_b)$ by Lemma~\ref{lem:handlenumberofboundarysumSeifertsurfaces}. One also observes that $g(S)=g(S_a)+g(S_b)$.
    
    Now suppose for each $i=a,b$ that $F_i$ is a Seifert surface for $K_i$ that realizes $g_{\MN}(K_i)$; hence $h(S_i) \geq h(F_i)$.
 Again, $F_i$ must be incompressible.   Then  $F= F_a \natural F_b$ is an incompressible Seifert surface for $K$ and again $h(F)=h(F_a)+h(F_b)$ by Lemma~\ref{lem:handlenumberofboundarysumSeifertsurfaces}.   However, since $S$ realizes $g_{\MN}(K)$, $h(S)=h(K) \leq h(F)$ and hence $h(S_a)+h(S_b) \leq h(F_a)+h(F_b)$.  But since $h(S_i) \geq h(F_i)$, we must have $h(S_i)=h(F_i)$.  Hence $g(F_i) \leq g(S_i)$ and so $g(F) \leq g(S)$.  Yet since $h(S_i)=h(F_i)$, Lemma~\ref{lem:handlenumberofboundarysumSeifertsurfaces} also gives that $h(S)=h(F)$.   Since $S$ realizes $g_{\MN}(K)$, we have $g(F) \geq g(S)$.  Now it follows that $g(F)=g(S)$ and $g(F_i)=g(S_i)$ as well.  Thus for each $i=a,b$ we have that $S_i$ realizes $g_{\MN}(K_i)$ too.
\end{proof}

\begin{theorem}\label{thm:additivityofGenusgMorseNovikov}
The minimal genus Morse-Novikov number is additive under connected sum.

    Given knots $K_a$ and $K_b$, 
    \[\MN_{g(K_a \# K_b)}(K_a \# K_b) = \MN_{g(K_a)}(K_a) + \MN_{g(K_b)}(K_b).\]
    
More specifically, if a Seifert surface $S$ for $K=K_a \# K_b$ realizes $\MN_{g(K)}(K)$, then $S = S_a \sharp S_b$ is the boundary sum of Seifert surfaces $S_i$ for $K_i$ that realize $\MN_{g(K_i)}(K_i)$ for each $i=a,b$. 
\end{theorem}

\begin{proof}
    Suppose $S$ is a minimal genus Seifert surface for $K = K_a \# K_b$  such that $h(S) = MN_{g(K)}(K)$.  By the additivity of Seifert genus under connected sum, $S$ decomposes as a boundary sum $S_a \sharp S_b$ for minimal genus Seifert surfaces $S_a$ and $S_b$ of the knots $K_a$ and $K_b$.  By Lemma~\ref{lem:handlenumberofboundarysumSeifertsurfaces}  we have 
    \[\MN_{g(K)}(K) = h(S) = h(S_a)+h(S_b) \geq \MN_{g(K_a)}(K_a) + \MN_{g(K_b)}(K_b).\]
    
    Now suppose for each $i=a,b$, let $S_i$ be a minimal genus Seifert surface for $K_i$ such that $h(S_i) = MN_{g(K_i)}(K_i)$.  Then $S_a \sharp S_b$ is a minimal genus Seifert surface for $K=K_a \# K_b$ by the additivity of Seifert genus.  Hence 
    \[\MN_{g(K)}(K) \leq h(S_a \sharp S_b) = h(S_a) + h(S_b) = \MN_{g(K_a)}(K_a) + \MN_{g(K_b)}(K_b)\]
    where the first equality is by Lemma~\ref{lem:handlenumberofboundarysumSeifertsurfaces}. 
\end{proof}

\begin{theorem}
    For each positive integer $n$, there exists a knot $K_n$ with 
    $g_\MN(K_n) - g(K_n) \geq n$ and
    $\MN_{g(K_n)}(K_n) - \MN(K_n) \geq 2n$.
\end{theorem}

\begin{proof}
    Let $K$ be a knot from Corollary~\ref{cor:Kangdoubles} where $g(K)=1$ while $g_\MN(K) \geq 2$.  For such knots we know that $\MN_{g(K)}(K)=\MN_1(K)=6$ and $2\leq \MN(K)\leq \MN_2(K)  \leq 4$.

    For any integer $n >0$, set $K_n = \#_n K$, the connected sum of $n$ copies of $K$.
    Then
    \[g(K_n) = n \cdot g(K) = n\]
    and 
    \[\MN_{g(K_n)}(K_n) = n \cdot \MN_{g(K)}(K) =6n\]
    Since Seifert genus and minimal genus Morse-Novikov numbers  are additive under connected sum, the latter by Theorem~\ref{thm:additivityofGenusgMorseNovikov}.
    On the other hand,
    \[g_\MN(K_n) = n \cdot g_\MN(K) \geq 2n\]
and
    \[\MN(K_n) = n \cdot \MN(K) \leq 4n\]
        since the Morse-Novikov genus and the Morse-Novikov number are additive under connected sum by Theorem~\ref{thm:MNgenusisadditive} and \cite{BakerAdditivityofMorseNovikov}.
Hence we have
\[g_\MN(K_n) - g(K) \geq 2n - n = n\]
        and
    \[\MN_{g(K_n)}(K_n) - \MN(K_n) \geq 6n - 4n = 2n\]
as desired.
\end{proof}

\begin{lemma}\label{lem:MNg-structure}
Let $L$ be an oriented link.
If $L$ is not fibered, then 
\begin{itemize}
    \item $\MN_{g} \geq \MN_{g+1}$ for all $g \geq g(L)$,
    \item $\MN_{g} = \MN$ for all $g \geq g_\MN(L)$. 
\end{itemize}
If $L$ is fibered, then $\MN_{g(L)} = \MN = 0$ and $\MN_{g} = 2$ for all $g > g(L)$.
\end{lemma}

\begin{proof}
    Given a circular Heegaard splitting $(M, F, S)$ with non-zero handle number, then we may create a new circular Heegaard splitting $(M,F_1', S'')$ of the same handle number where $g(F_1')= g(F)+1$ 
    using the tubing operation of \cite[Section 3.4]{BMG-MorseNovikovHandleNumbers} as follows.  Trivially ``inflate'' $(M,F,S)$ to a generalized circular Heegaard splitting $(M, F_0 \cup F_1, S_0 \cup S_1)$ where $F$ is thickened to a product $F \times [0,1]$, $F_0 = F\times \{0\}$, $S_0 = F \times \{\frac12\}$, $F_1 = F \times \{1\}$, and $S=S_1$.   Then a $1$-handle of $(M,F,S)$ may be identified with a $1$-handle in the compression body $A_1$ between $F_1$ and $S_1$.  Apply the tubing operation of \cite[Section 3.4]{BMG-MorseNovikovHandleNumbers} where $i=1$ and $j=k=0$ to obtain a generalized circular Heegaard splitting $(M, F_0 \cup F'_1, S'_0 \cup S_1)$.  Finally, amalgamate along $F_0$ to obtain $(M, F'_1, S'')$.  By \cite[Lemmas 3.1 and 3.7]{BMG-MorseNovikovHandleNumbers}, $h(M, F, S) \geq h(M, F_1', S'')$. By construction $F_1'$ is just a tubing of $F$, so $g(F_1')= g(F)+1$. 

    Consequently, if a genus $g$ Seifert surface $F$ realizes $\MN_{g}$ so that $\MN_g= h(F)$, then a tubing as above creates a Seifert surface $F_1'$ of genus $g+1$ with $h(F) \geq h(F_1')$.  Hence $\MN_g$ is a non-decreasing function of $g$ down to $\MN$ except when $\MN_g=0$ in which case there are no $1$-handles for the tubing operation.
    
    Indeed, such an exception only occurs when the link is fibered and $g$ is the Seifert genus.  Then a trivial stabilization of a fiber has genus $g+1$ and handle number $2$.  From there, all tubings must have handle number $2$ as well since the fiber is the unique incompressible Seifert surface.
\end{proof}

\section{Problems}

All the knots we have constructed here with Morse-Novikov genus greater than Seifert genus are satellited knots, namely twisted Whitehead doubles and connected sums thereof. Certainly the same phenomenon occurs for hyperbolic knots. (It cannot for torus knots since torus knots are all fibered.)

\begin{prob}\label{prob:hypexamples}
    Show there exists a hyperbolic knot $K$ with $g_\MN(K) > g(K)$.  Indeed, show the knots $K_n$ of Theorem~\ref{thm:arblargediscrepancy} can be chosen to all be hyperbolic.
\end{prob}  

As suggested in the introduction, for the known examples of knots with Morse-Novikov number $2$, the Seifert genus is the Morse-Novikov genus. Is this always the case?  Indeed, as stated in Corollary~\ref{cor:Kangdoubles}, we only show the Morse-Novikov number is either $2$ or $4$ for our examples constructed here. 

\begin{prob}
    If a knot has Morse-Novikov number $2$, must the Seifert genus equal the Morse-Novikov genus?
\end{prob}

Lemma~\ref{lem:MNg-structure} shows that for any given non-fibered knot,  $\MN_g$ is non-increasing as a function of $g$.  For our knots $K$ constructed here, we have shown that  $\MN_g$ takes at least two values: $\MN_{g(K)}$ and $\MN_{g_\MN(K)}$. Surely there are knots for which $\MN_g$ takes on more than two values.

\begin{prob}
    Find a sequence of knots $K_n$ so that $\MN_g(K_n)$, as a function of $g$, takes on at least $n$ values.
\end{prob} 

Furthermore, what can we say about where these drops in value for $\MN_g$ occur?  
\begin{prob}\label{prob:decreaseatincomp}
    If $\MN_g > \MN_{g+1}$ then is $\MN_{g+1}$ realized by an incompressible Seifert surface?
\end{prob}

It may be more fruitful to work counter to this.
The tubing operation used in the proof of Lemma~\ref{lem:MNg-structure} suggests one avenue for constructing a negative answer to the above Problem~\ref{prob:decreaseatincomp}.

\begin{prob}
    Find a Seifert surface $F$ that compresses to a Seifert surface $F'$ with larger handle number, ie.\ with $h(F) < h(F')$.
\end{prob}

\bibliographystyle{alpha}
\bibliography{biblio}

\begin{thebibliography}{MGNnRL15}

\bibitem[Bak21]{BakerAdditivityofMorseNovikov}
Kenneth~L. Baker.
\newblock The {M}orse-{N}ovikov number of knots under connected sum and cabling.
\newblock {\em J. Topol.}, 14(4):1351--1368, 2021.

\bibitem[BMG22]{BMG-MorseNovikovHandleNumbers}
Kenneth~L. Baker and Fabiola Manjarrez-Gutiérrez.
\newblock Morse-novikov numbers, tunnel numbers, and handle numbers of sutured manifolds, 2022.

\bibitem[BMG23]{BMG-nearlyfibered}
Kenneth~L. Baker and Fabiola Manjarrez-Gutiérrez.
\newblock Handle numbers of guts of sutured manifolds and nearly fibered knots, 2023.

\bibitem[God92]{goda-murasugisum}
Hiroshi Goda.
\newblock Heegaard splitting for sutured manifolds and {M}urasugi sum.
\newblock {\em Osaka J. Math.}, 29(1):21--40, 1992.

\bibitem[God93]{goda-handlenumber}
Hiroshi Goda.
\newblock On handle number of {S}eifert surfaces in {$S^3$}.
\newblock {\em Osaka J. Math.}, 30(1):63--80, 1993.

\bibitem[Kak92]{KakimizuComplex}
Osamu Kakimizu.
\newblock Finding disjoint incompressible spanning surfaces for a link.
\newblock {\em Hiroshima Math. J.}, 22(2):225--236, 1992.

\bibitem[Kan15]{Kang}
Sungmo Kang.
\newblock Knots admitting {S}eifert-fibered surgeries over {$S^2$} with four exceptional fibers.
\newblock {\em Bull. Korean Math. Soc.}, 52(1):313--321, 2015.

\bibitem[Lyo74]{lyon}
Herbert~C. Lyon.
\newblock Simple knots without unique minimal surfaces.
\newblock {\em Proc. Amer. Math. Soc.}, 43:449--454, 1974.

\bibitem[MGNnRL15]{freegenusone}
Fabiola Manjarrez-Guti\'errez, V\'ictor N\'u\~nez, and Enrique Ram\'irez-Losada.
\newblock Circular handle decompositions of free genus one knots.
\newblock {\em Pacific J. Math.}, 275(2):361--407, 2015.

\bibitem[Mor89]{Morimoto}
Kanji Morimoto.
\newblock On minimum genus {H}eegaard splittings of some orientable closed {$3$}-manifolds.
\newblock {\em Tokyo J. Math.}, 12(2):321--355, 1989.

\end{thebibliography}

\end{document}